\documentclass[reqno]{amsart}

\usepackage{enumerate}
\usepackage{amsmath}%
\usepackage{amsfonts}%
\usepackage{amssymb}%
\usepackage{graphicx}
\usepackage{mathrsfs}
\usepackage{hyperref}
\usepackage{fullpage}
\usepackage{lineno}
\usepackage{setspace}
\usepackage{enumitem}
\newtheorem{theorem}{Theorem}
\theoremstyle{plain}

\newtheorem{construction}[theorem]{Construction}

\newtheorem{lemma}[theorem]{Lemma}

\newtheorem{proposition}[theorem]{Proposition}

\numberwithin{equation}{section}
\numberwithin{theorem}{section}
\numberwithin{case}{section}

\numberwithin{subcase}{case}

\def\F{\mathcal{F}}

\def \a{\alpha}
\def \e{\epsilon}
\def \r{\gamma}
\def \cP{\mathcal{P}}

\def \bfi{\mathbf{i}}
\def \bfu{\mathbf{u}}
\def \bfv{\mathbf{v}}
\def \bfw{\mathbf{w}}

\begin{document}
\title{On Perfect Matchings in $k$-complexes}
\author{Jie Han}
\thanks{Research partially supported by Simons Foundation \#630884.}

\address{Department of Mathematics, University of Rhode Island, 5 Lippitt Road, Kingston, RI, 02881, USA}
\email{jie\_han@uri.edu}
\maketitle

\onehalfspace

\begin{abstract}
Keevash and Mycroft [\emph{Mem.~Amer.~Math.~Soc., 2015}] developed a geometric theory for hypergraph matchings and characterized the dense simplicial complexes that contain a perfect matching.
Their proof uses the hypergraph regularity method and the hypergraph blow-up lemma recently developed by Keevash. 
In this note we give a new proof of their results, which avoids these complex tools.
In particular, our proof uses the lattice-based absorbing method developed by the author and a recent probabilistic argument of Kohayakawa, Person and the author.
\end{abstract}

\section{Introduction}

Matchings are fundamental objects in graph theory and have broad practical applications in other branches of science and a variety of real-world problems.
Perhaps the best known application of these is in the assignment of graduating medical students to their first hospital appointments\footnote{In 2012, the Nobel Memorial Prize in Economics was awarded to Shapley and Roth ``for the theory of stable allocations and the practice of market design."}.
Matchings in hypergraphs also offer a universal framework for many important combinatorial problems, e.g., the existence conjecture for block designs\footnote{This conjecture was recently solved by Keevash~\cite{Keevash_design}. See~\cite{GKLO} for another proof given by Glock et al.} and Ryser's conjecture on transversals in Latin square.
Moreover, there are examples for applications to practical problems, e.g., the `Santa Claus' allocation problem~\cite{santa_claus}.

In this paper, we are particularly interested in the existence of \emph{perfect matchings}.
A matching is perfect if it covers all vertices of the host graph.
Perfect matchings in graphs are well understood.
For instance, Tutte's Theorem \cite{Tu47} gives necessary and sufficient conditions for a graph to contain a perfect matching, and Edmonds' Algorithm \cite{Edmonds} finds such a matching in polynomial time. 
However, as probably the most straightforward objects one can extend to hypergraphs, perfect matchings in $k$-uniform hypergraphs appear to be considerably harder than the graph case.
In fact, for $k\ge 3$, determining whether a $k$-uniform hypergraph contains a perfect matching was one of Karp's celebrated 21 NP-complete problems \cite{Karp} and by now no characterization theorem (such as Tutte's Theorem in the graph case) is known.

Since the general problem is intractable provided P $\neq$ NP, it is natural to ask for conditions on $H$ which guarantee the existence of a perfect matching. One well-studied class of such conditions are minimum degree conditions. 
Recently, there has been a strong interest on matching problems in dense hypergraphs, due to the recent development in techniques such as the hypergraph regularity method and the absorption method.

A \emph{hypergraph} $H$ consists of a vertex set $V(H)$ and an edge set $E(H)\subseteq 2^{V(H)}$.
Given $k\ge 2$, a \emph{$k$-uniform hypergraph} (or $k$-graph) $H$ is a hypergraph whose every edge has cardinality $k$.
In this note we often identify $H$ and its edge set $E(H)$. 
A \emph{matching} in $H$ is a collection of vertex-disjoint edges of $H$. A \emph{perfect matching} $M$ in $H$ is a matching that covers all vertices of $H$. 
The \emph{minimum codegree} of a $k$-uniform hypergraph $H$ is the maximum integer $t$ such that every $(k-1)$-vertex subset of $H$ is in at least $t$ edges.
R\"odl, Ruci\'nski and Szemer\'edi \cite{RRS09} determined the sharp minimum codegree condition that ensures a perfect matching in a $k$-uniform hypergraph on $n$ vertices for large $n$ and all $k\ge 3$. 
The value is ${n}/{2}-k+C$, where $C\in\{3/2, 2, 5/2, 3\}$ depends on the values of $n$ and $k$. 
Since then there have been a large number of efforts from various researchers on the minimum degree conditions that force perfect matchings in $k$-uniform hypergraphs, see~\cite{AFHRRS, CzKa, HPS, Khan2, Khan1, KO06mat, KOT, MaRu, Pik, RRS06mat, RRS09, TrZh12, TrZh13, TrZh15} and the recent surveys \cite{RR, zsurvey}. 

In this paper we consider a more general setting instead of uniform hypergraphs.
For $k\ge 2$, a \emph{$k$-complex} is a hypergraph $J$ such that every edge $e$ of $J$ contains at most $k$ vertices and every subset of $e$ is also an edge of $J$.
%
A celebrated result \cite{KM1} of Keevash and Mycroft characterizes the (reasonably) dense $k$-complexes that contain perfect matchings in its $k$-th level.
It can be viewed as an analogue of Tutte's Theorem into dense hypergraphs.
Since the statement of the characterizations are somewhat technical, we defer them to Section 2.
Their result also has exciting applications.
Here is a non-exhaustive list:
\begin{itemize}
\item Using their characterization theorem, Keevash and Mycroft \cite{KM1} determined the exact codegree threshold for tiling tetrahedron\footnote{The tetrahedron is the unique $3$-graph with 4 vertices and 4 edges; an \emph{$F$-tiling} is a spanning subgraph consisting of vertex-disjoint copies of $F$. \label{ft:1}} in 3-uniform hypergraphs.
\item In a subsequent paper \cite{KeMy2}, they used the multi-partite version of the characterization theorem to establish a multi-partite analogue of the celebrated Hajnal--Szemer\'edi Theorem.
\item Moreover, Keevash, Knox and Mycroft~\cite{KKM13} solved a hardness problem for perfect matchings in $k$-uniform hypergraphs posed by Karpi\'nski, Ruci\'nski and Szyma\'nska~\cite{KRS10} almost entirely (the author~\cite{Han14_poly} recently solved the only missing case).
\item Han et~al.~\cite{HLTZ_K4} used a result on almost perfect matching~\cite[Lemma 7.6]{KM1} (see Lemma~\ref{lem:keevash} below) to determine the sharp minimum codegree condition on tiling $K_4^{3-}$. \footnote{A $K_4^{3-}$ is the unique $3$-graph with 4 vertices and 3 edges.}
\end{itemize}

The proof in~\cite{KM1} used the \emph{hypergraph regularity method} and the \emph{blow-up lemma for hypergraphs} developed recently by Keevash.
The hypergraph regularity method -- the extension of Szemer\'edi's graph regularity lemma~\cite{Sze} to the setting of $k$-graphs -- is one of the most celebrated combinatorial results in this century. 
By now there are several (very different) proofs of this lemma, obtained by Gowers~\cite{Gowers}, by Nagle--R\"odl--Schacht--Skokan~\cite{NRS06,RoSk,RS, RoSc_regularity} and by Tao~\cite{Tao}.
The blow-up lemma, developed by~Koml\'os, Sark\"ozy and Szemer\'edi~\cite{Blowup}, is a celebrated tool for using the regularity lemma to embed spanning subgraphs of bounded degree.
It has many applications to embedding spanning subgraphs, including the famous Pos\'a--Seymour conjecture~\cite{KSS-posa}, see the survey~\cite{KuOs-survey}.
However, there are two main drawbacks of the regularity--blowup method.
Firstly, the employment of the regularity--blowup method usually makes the arguments quite involved and technical; secondly, the regularity method requires the order of the host graph to be extremely large -- e.g., when it is applied to $k$-graphs, the order is given by the $k$-th Ackermann function.

The goal of this note is to give new and short proofs of the characterization theorems in~\cite{KM1} (see Theorems~\ref{thm:main},~\ref{thm:main2} and~\ref{thm:711} below).
The absorbing method, initiated by R\"odl, Ruci\'nski and Szemer\'edi \cite{RRS06}, has proven to be efficient in embedding spanning structures in graphs and hypergraphs. 
More precisely, (in case for perfect matchings) it reduces the problem into two subproblems:
\begin{itemize}
\item[(1)] find an almost perfect matching, that is, a matching covering all but a $o(1)$ proportion of vertices;
\item[(2)] turn the almost perfect matching to a perfect matching.
\end{itemize}
To achieve (2) we will use a novel variant, namely, the \emph{lattice-based absorbing method} developed recently by the author \cite{Han14_poly}.
We note that the degree assumptions in our problem are too low to apply the standard absorbing method. 
Secondly, to achieve (1), we will use a recent probabilistic approach developed by Kohayakawa, Person and the author~\cite{HKP} to `convert' a class of `weight-disjoint' perfect fractional matchings to an almost perfect matching.
To find a perfect fractional matching, we use a result~\cite[Theorem 7.2]{KM1}, which is a careful application of the well-known Farkas' Lemma.

Our entire proof is regularity-free, and also avoid the use of the hypergraph blow-up lemma.
This allows us to reduce the technicality in the proofs greatly, compared with the original proof in~\cite{KM1}.
Secondly, our proof mainly use probabilistic methods, which usually yield moderate bounds on the order of the host graphs (see Concluding Remarks).
At last, using our result, the proofs given in~\cite{KM1, KeMy2, HLTZ_K4} are now regularity-free\footnote{The proof of~\cite{KKM13} still uses the weak regularity lemma as a tool, which is considered as a much lighter tool compared with the (strong) regularity lemma used in~\cite{KM1}.}.

\section{The characterization theorems}
Let $J$ be a $k$-complex.
For $0\le i\le k$, let $J_i\subseteq J$ be the set of edges of size exactly $i$. 
Then each $J_i$ is $i$-uniform and $J=\bigcup_{i=0}^k J_i$.
Fix an $i$-edge $e\in J_i$, the \emph{degree $\deg_J(e)$} is the number of $(i+1)$-edges $e'$ of $J_{i+1}$ which contain $e$ as a subset (note that this is \emph{not} the standard notion of degree used in $k$-graphs, in which the degree of a set is the number of edges containing it).
For $0\le i<k$, the \emph{minimum $i$-degree} of $J$, denoted by $\delta_i(J)$, is the minimum $\deg_J(e)$ taken over all $i$-edges $e$ of $J_i$.
The \emph{degree sequence of $J$} is $\delta(J)=(\delta_0(J), \delta_1(J), \dots, \delta_{k-1}(J))$. 

\subsection{The barriers that prevent the existence of a perfect matching}

We start with some very natural constructions that prevent the existence of a perfect matching in dense hypergraphs/complexes, which are important to us.
\begin{construction}[Space Barrier, \cite{KM1}]\label{con:sb}
Let $V$ be a set of size $n$, $j\in [k-1]$ and $S\subseteq V$. Let $J=J(S, j)$ be the $k$-complex in which for every $i\in[k]$, $J_i$ consists of all $i$-sets in $V$ that contain at most $j$ vertices of $S$.
Since each $k$-edge contains at most $j$ vertices of $S$, if $|S|> jn/k$, then $J_k$ contains no perfect matching.
\end{construction}

For $j\in [k-1]$, the degree sequence of $J=J(S, j)$ is
\begin{equation}
\delta(J)= \left( n, n-1,\dots, n-(j-1), n-|S|, n-|S|-1, \dots, n-|S|-(k-j-1)  \right). \nonumber
\end{equation}
Thus, to force the existence of a perfect matching, this suggests a degree sequence `at least'
\begin{equation}
\delta(J)\ge \left( n, \frac{k-1}{k} n, \frac {k-2}k n, \dots, \frac 1k n  \right), \label{eq:sb}
\end{equation}
that is, each individual digit in the degree sequence above cannot be lowered.
However, we will see below that another class of barriers of considerably higher degrees also prevent the presence of a perfect matching.

We also review the divisibility barriers that are observed in \cite{TrZh12} and then generalized in \cite{KM1}. Let $V$ be a set of vertices, and let $\mathcal P$ be a partition of $V$ into $r$ parts $V_1,\dots, V_r$. 
The \emph{index vector} $\mathbf{i}_{\cP}(S)\in \mathbb{Z}^r$ of a subset $S\subseteq V(H)$ with respect to $\cP$ is the vector whose coordinates are the sizes of the intersections of $S$ with each part of $\cP$, i.e., $\mathbf{i}_{\cP}(S)_{V_i}=|S\cap V_i|$ for $i\in [r]$.
Throughout this note, every partition has an implicit order on its parts. 

\begin{construction}[Divisibility barrier, \cite{TrZh12, KM1}]
Let $\cP$ partition a vertex set $V$ into $d$ parts.
Suppose $L$ is a lattice in $\mathbb Z^d$ with $\bfi_\cP(V)\notin L$.
Fix any $k\ge 2$, and let $H$ be the $k$-graph on $V$ whose edges are all $k$-tuples $e$ with $\bfi_\cP(e)\in L$.
For any matching $M$ in $H$ with vertex set $S=\bigcup_{e\in M}e$ we have $\bfi_\cP(S)=\sum_{e\in M} \bfi_\cP(e)\in L$. Since $\bfi_\cP(V)\notin L$ it follows that $H$ does not have a perfect matching.
\end{construction}

The following example is a special case of the divisibility barriers that have been observed earlier.
Let $H$ be a 3-graph defined as follows. 
Let $A\cup B$ be a partition of $V(H)$ and let the edge set of $H$ consist of all triples that intersecting $B$ at an even number of vertices.
If $|B|$ is odd, then $H$ has no perfect matching.
Indeed, one can realize this as the divisibility barrier by letting $\mathcal P=(A, B)$, $L=\langle (1,2), (3,0) \rangle\subseteq \mathbb Z^2$ and $\bfi_\cP(V)\notin L$.
At last, note that the degree sequence of the 3-complex induced by $H$ is at least $(n, n-1, n/2-1)$.

\subsection{The characterization theorems}
In fact, Keevash and Mycroft provided a series of characterization theorems in~\cite[Section 2]{KM1} and below we only state two main ones.
However, we will reprove their main technical result, which does recover all of the results in~\cite[Section 2]{KM1} via a regularity-free approach.
Roughly speaking, the main result in \cite{KM1} says that if a dense $k$-complex $J$ is not `close' to either the space barriers or the divisibility barriers, then $J_k$ has a perfect matching.
To describe the `closeness', we use the following definitions from \cite{KM1}.

Fix an integer $r>0$, let $H$ be a $k$-graph and let $\cP=\{V_1,\dots, V_r\}$ be a partition of $V(H)$.
We call a vector $\mathbf{i}\in \mathbb{Z}^r$ an \emph{$s$-vector} if all its coordinates are nonnegative and their sum equals $s$.
Given $\mu>0$, a $k$-vector $\mathbf{v}$ is called a $\mu$\emph{-robust edge-vector} if at least $\mu |V(H)|^k$ edges $e\in E(H)$ satisfy $\mathbf{i}_\cP(e)=\mathbf{v}$.
Let $I_{\cP}^{\mu}(H)$ be the set of all $\mu$-robust edge-vectors and let $L_{\cP}^{\mu}(H)$ be the lattice (additive subgroup) generated by the vectors of $I_{\cP}^{\mu}(H)$.
A lattice $L\subseteq \mathbb Z^r$ is called \emph{complete} (otherwise, \emph{incomplete}) if it contains all $k$-vectors in $\mathbb Z^r$.
For $j\in [r]$, let $\mathbf{u}_j\in \mathbb{Z}^r$ be the $j$-th \emph{unit vector}, namely, $\mathbf{u}_j$ has 1 on the $j$-th coordinate and 0 on other coordinates.
A \emph{transferral} is the vector $\bfu_i - \bfu_j$ for some $i\neq j$.

Now we are ready to state one of their main theorems~\cite[Theorem 2.9]{KM1}.
Throughout this paper, $x\ll y$ means that for any $y\ge 0$ there exists $x_0> 0$ such that for any $0<x\le x_0$ the subsequent statement holds.
Hierarchies of other lengths are defined similarly.

\begin{theorem}\cite[Theorem 2.9]{KM1}\label{thm:main}
Suppose that $1/n\ll \r \ll \mu, \beta \ll 1/k$ and that $k$ divides $n$. Let $J$ be a $k$-complex on $n$ vertices such that
\begin{equation}\label{eq:deg}
\delta(J)\ge \left( n, \left(\frac{k-1}{k} - \r\right) n, \left(\frac {k-2}k - \r\right) n, \dots, \left(\frac 1k -\r \right) n  \right).
\end{equation}
Then $J$ has at least one of the following properties:
\begin{enumerate}[label=$(\roman*)$]
\item {\bf Matching:} $J_k$ contains a perfect matching.
\item {\bf Space barrier:} For some $p \in [k-1]$ and set $S \subseteq V$ with $|S| = \lfloor pn/k \rfloor$, we have $e(J_{p+1}[S] )\le \beta n^{p+1}$.
\item {\bf Divisibility barrier:} There is some partition $\cP$ of $V(J)$ into $d\le k$ parts of size at least $\delta_{k-1}(J) - \mu n$ such that $L_{\cP}^{\mu}(J_k)$ is incomplete and transferral-free.
\end{enumerate}
\end{theorem}

Note that the degree sequence condition~\eqref{eq:deg} is slightly weaker than that given by the space barriers as in~\eqref{eq:sb}.

The second theorem is a multi-partite analogue of Theorem~\ref{thm:main}. To state the theorem, we need some more definitions from~\cite{KM1}.
Let $H$ be a hypergraph, and let $\cP$ be a partition of $V(H)$ into $V_1,\dots, V_r$.
Then we say a set $S$ of vertices and its index vector $\bfi_\cP(S)$ are \emph{$\cP$-partite} if $S$ has at most one vertex in any part of $\cP$.
We say $H$ is \emph{$\cP$-partite} if every edge of $H$ is $\cP$-partite.
We use $\binom{\cP}k$ to denote the set of all $\cP$-partite $k$-vectors.
Given a partite $k$-complex, we define the following alternative notion of degree.
Let $V$ be a set of vertices, let $\cP$ be a partition of $V$ into $r$ parts $V_1,\dots, V_r$, and let $J$ be a $\cP$-partite $k$-complex on $V$.
For each $0\le j\le k-1$ we define the \emph{partite minimum $j$-degree $\delta_j^*(J)$} as the largest $m$ such that any $j$-edge $e$ has at least $m$ extensions to a $(j+1)$-edge in any part not used by $e$, that is,
\[
\delta_j^*(J) := \min_{e\in J_j} \min_{i: e\cap V_i = \emptyset} |\{v\in V_i: e\cup \{v\}\in J\}|.
\]
The \emph{partite degree sequence} is $\delta^*(J) = (\delta_0^*(J),\dots, \delta_{k-1}^*(J))$.

For a matching $M$ in a $\cP$-partite $k$-graph $H$ we write $n_{\bfi}(M)$ to denote the number of edges in $M$ with index vector $\bfi$.
We say that $M$ is balanced if $n_\bfi(M)$ is constant over all $\cP$-partite $k$-vectors $\bfi$.
However, \cite[Construction 2.11]{KM1} shows that one cannot guarantee a balanced matching in the partite analogue of Theorem~\ref{thm:main}.
Instead a weaker property was introduced in \cite{KM1}: we say that $M$ is \emph{$\a$-balanced} if $n_\bfi(M) \ge (1-\a)n_{\bfi'}(M)$ for any two $\cP$-partite $k$-vectors $\bfi, \bfi'\in \binom{\cP}k$.

At last, for two partitions $\cP, \cP'$ on the same vertex set, we say that $\cP'$ refines $\cP$ if every part of $\cP'$ is a subset of some part of $\cP$.
Given two partitions $\cP, \cP'$ such that $\cP'$ refines $\cP$, a lattice $L$ on $\cP'$ is \emph{complete (otherwise incomplete) with respect to $\cP$} if it contains all the $\cP'$-partite $k$-vectors; a lattice $L$ on $\cP'$ is \emph{transferral-free with respect to $\cP$} if $L$ does not contain a transferral $\bfu_i - \bfu_j$ such that $V_i \cup V_j$ is a subset of some part of $\cP$.

\begin{theorem}\cite[Theorem 2.10]{KM1}\label{thm:main2}
Suppose that $1/n\ll \r, \a \ll \mu, \beta \ll 1/r \le 1/k$.
Let $\cP$ partition a set $V$ into parts $V_1,\dots, V_r$ each of size $n$, where $k\mid r n$.
Suppose that $J$ is a $\cP$-partite $k$-complex with
\begin{equation*}
\delta^*(J)\ge \left( n, \left(\frac{k-1}{k} - \r\right) n, \left(\frac {k-2}k - \r\right) n, \dots, \left(\frac 1k -\r \right) n  \right).
\end{equation*}
Then $J$ has at least one of the following properties:
\begin{enumerate}[label=$(\roman*)$]
\item {\bf Matching:} $J_k$ contains an $\a$-balanced perfect matching.
\item {\bf Space barrier:} For some $p \in [k-1]$ and sets $S_i \subseteq V_i$ with $S_i = \lfloor pn/k \rfloor$ for all $i\in [r]$, we have $e(J_{p+1}[S] )\le \beta n^{p+1}$, where $S:=\bigcup_{i\in [r]}S_i$.
\item {\bf Divisibility barrier:} There is some partition $\cP'$ of $V(J)$ into $d\le k r$ parts of size at least $\delta^*_{k-1}(J) - \mu n$ such that $\cP'$ refines $\cP$ and $L_{\cP'}^{\mu}(J_k)$ is incomplete and transferral-free with respect to $\cP$.
\end{enumerate}
\end{theorem}



Instead of proving Theorems~\ref{thm:main} and~\ref{thm:main2} directly, we will (re)prove the main technical result, namely, \cite[Theorem 7.11]{KM1}, from which Theorems~\ref{thm:main} and~\ref{thm:main2} can be simply derived as in~\cite{KM1}.

\section{The main technical result}

To state the technical result of~\cite{KM1}, we need some further definitions from~\cite{KM1}.
Let $k$ and $r$ be positive integers.
An allocation function $f$ is a function $f:[k]\to [r]$.
The index vector of $f$ is defined as $\bfi(f):=(|f^{-1}(1)|, \dots, |f^{-1}(r)|)\in \mathbb Z^r$.

Let $\cP$ be a partition of $r$ parts and $I$ be a multiset of $k$-vectors.
Then we may form a multiset $F$ of allocation functions $f:[k]\to [r]$ as follows.
For each $\bfi\in I$ (with repetition) choose an allocation $f$ with $\bfi(f)=\bfi$, and include in $F$ each of the $k!$ permutations $f^\sigma$ for $\sigma\in Sym_k$ (again including repetitions), where $f^\sigma(i) = f(\sigma(i))$ for $i\in [k]$.
Note that the multiset $F$ so obtained does not depend on the choices of allocation function $f$.
If $F$ can be obtained this way, we call $F$ an \emph{allocation} and we write $I(F)$ for the multiset $I$ from which $F$ was defined.
We say that an allocation $F$ is \emph{$(k,r)$-uniform} if for every $i\in [k]$ and $j\in [r]$ there are $|F|/r$ functions $f\in F$ with $f(i)=j$.
We also say that $F$ is \emph{connected} if there is a connected graph $G_F$ on $[r]$ such that for every $jj'\in E(G_F)$ and $i, i'\in [k]$ with $i \neq i'$ there is $f\in F$ with $f(i) = j$ and $f(i') = j'$.
We remark that although clearly there are at most $r^k$ distinct elements in $F$, as a multiset (note that $I$ is a multiset as well) there is no natural upper bound on $|F|$.
So in the result we require an upper bound as $|F|\le D_F$ for some constant $D_F>0$.

A \emph{$k$-system} is a hypergraph with all edges of size at most $k$; namely, it is not necessarily 'downward-closed' as $k$-complexes.
For an allocation $F$, a $k$-system $J$ on $V$ is $\cP F$-partite if for any $j\in [k]$ and $e\in J_j$ there is some $f\in F$ so that $e = \{v_1, \dots, v_j\}$ with $v_i\in V_{f(i)}$ for $i\in [j]$, namely, every edge of $J$ can be constructed through the process above for some $f\in F$.
The minimum $F$-degree sequence of $J$ is then defined to be
\[
\delta^F(J) := (\delta_0^F(J), \dots, \delta_{k-1}^F(J)),
\]
where $\delta_j^F(J):=\min_{f\in F}\delta_j^f(J)$ and $\delta_j^f(J)$ is the largest $m$ such that for any $\{v_1,\dots, v_j\}\in E(J)$ with $v_i\in V_{f(i)}$ for $i\in [j]$ there are at least $m$ vertices $v_{j+1}\in V_{f(j+1)}$ such that $\{v_1,\dots, v_{j+1}\}\in E(J)$.
Note that the minimum $F$-degree sequence generalizes simultaneously the minimum degree sequence ($r=1$) and the $r$-partite minimum degree sequence ($r\ge k$ and $F$ be the collection of all injections from $[k]$ to $[r]$).

Our next definition generalizes the notion of `$\a$-balancedness'.
Let $J$ be a $k$-system, $M$ be a perfect matching of $J_k$ and $F$ be a $(k,r)$-uniform allocation.
We say that $M$ \emph{$\alpha$-represents $F$} if for any $\bfi, \bfi'\in I(F)$, letting $n_{\bfi}(M)$, $n_{\bfi'}(M)$ denote the number of edges $e'\in M$ with index $\bfi, \bfi'$ divided by the multiplicities of $\bfi, \bfi'$ respectively\footnote{Note that different allocation functions $f$ in $F$ may correspond to the same index vector, so when we consider the `balancedness' we have to divide the number of edges of index vector $\bfi$ by the multiplicity of $\bfi$ in $F$.}, we have $n_{\bfi'}(M)\ge (1-\alpha)n_{\bfi}(M)$.
When $M$ 0-represents $F$, we say $M$ is \emph{$F$-balanced}.

Now we are ready to state~\cite[Theorem 7.11]{KM1}.
In the following, we assume the following degree sequence condition (for a certain choice of $F$):
\begin{equation}\label{eq:deg3}
\delta^F(J)\ge \left( n, \left(\frac{k-1}{k} - \r\right) n, \left(\frac {k-2}k - \r\right) n, \dots, \left(\frac 1k -\r \right) n  \right).
\end{equation}

\begin{theorem}\label{thm:711}\cite{KM1}
Let $1/n \ll \alpha \ll \gamma \ll\beta , \mu \ll 1/D_F, 1/r, 1/k$.
Suppose $F$ is a $(k,r)$-uniform connected allocation with $|F|\le D_F$, and that $k$ divides $r n$.
Let $\cP$ be a partition of a set $V$ into parts $V_1,\dots, V_r$ of size $n$, and $J$ be a $\cP F$-partite $k$-complex on $V$ satisfying~\eqref{eq:deg3} and
\begin{enumerate}[label=$(\roman*)$]
\item for any $p \in [k-1]$ and sets $S_i \subseteq V_i$ with $S_i = \lfloor pn/k \rfloor$ for all $i\in [r]$, we have $e(J_{p+1}[S] )> \beta n^{p+1}$, where $S:=\bigcup_{i\in [r]}S_i$. \footnote{In the original statement of \cite[Theorem 7.11]{KM1}, this item is stated as: for such $p$ and $S$ there are at least $\beta n^k$ edges of $J_k$ with more than $p$ vertices in $S$. Because of~\eqref{eq:deg3}, these two statements are equivalent up to a constant factor. Indeed, if $e(J_{p+1}[S] )> \beta n^{p+1}$, then one can grow these $(p+1)$-edges in $S$ to $k$-edges by the minimum $F$-degree condition greedily; for the converse, one can obtain $(p+1)$-edges in $S$ by averaging (e.g., dividing by $r^k n^{k-p-1}$).} \label{item:i}
\item $L_{\cP'}^{\mu}(J_k)$ is complete with respect to $\cP$ for any partition $\cP'$ of $V(J)$ which refines $\cP$ and whose parts each have size at least $\delta_{k-1}^F(J)-\mu n$. \label{item:iii}
\end{enumerate}
Then $J_k$ contains a perfect matching which $\alpha$-represents $F$.
\end{theorem}

Theorems~\ref{thm:main} and~\ref{thm:main2} follow from Theorem~\ref{thm:711} as shown in~\cite{KM1}.
We omit the deductions and therefore the rest of this note is devoted to the new proof of Theorem~\ref{thm:711}.

We use the following lemma proved by Keevash and Mycroft~\cite[Lemma 7.6]{KM1}.
Their proof uses the fractional matchings together with hypergraph regularity, which is a known way of `turning' a perfect fractional matching into an almost perfect matching.
In Section 5 we will give another proof of this result, which uses fractional matchings but together with a new probabilistic approach developed recently by Kohayakawa, Person and the author~\cite{HKP}.
Note that the following lemma was proved for \emph{$k$-systems}, namely, one does not need the downward-closed property of $k$-complexes.

\begin{lemma}\cite[Lemma 7.6]{KM1} \label{lem:keevash}
Suppose that $1/n \ll \phi \ll \r \ll \beta, 1/D_F, 1/r, 1/k$.
Suppose $F$ is a $(k,r)$-uniform connected allocation with $|F|\le D_F$.
Let $\cP$ be a partition of a set $V$ into parts $V_1,\dots, V_r$ of size $n$, and $J$ be a $\cP F$-partite $k$-system on $V$ 
satisfying \eqref{eq:deg3} and~\ref{item:i} in Theorem~\ref{thm:711}.
Then $J_k$ contains an $F$-balanced matching $M$ which covers all but at most $\phi n $ vertices of $V$.
\end{lemma}

Next we state our absorbing lemma which, in fact, works under an arbitrarily small minimum $F$-degree sequence condition.

\begin{lemma}[Absorbing Lemma]\label{lem:abs}
Suppose that $1/n \ll \phi \ll \e \ll \mu \ll \zeta, 1/D_F, 1/r, 1/k$. 
Suppose $F$ is a $(k,r)$-uniform connected allocation with $|F|\le D_F$.
Let $\cP$ be a partition of a set $V$ into parts $V_1,\dots, V_r$ of size $n$, and $J$ be a $\cP F$-partite $k$-complex on $V$ satisfying $\delta^F(J)\ge (n, \zeta n, \dots, \zeta n)$ and \ref{item:iii} in Theorem~\ref{thm:711}.
Then there exists a balanced set $W\subseteq V$ of order at most $\e n$ such that the following holds. 
Let $U \subseteq V\setminus W$ be any set such that $|U|\le \phi n$ and $k\mid |U|$.
Then both $J_k[W]$ and $J_k[U\cup W]$ contain perfect matchings.
\end{lemma}

We combine these two lemmas and give a new proof of Theorem~\ref{thm:711}.

\begin{proof}[Proofs of Theorem~\ref{thm:711}]
We apply Lemma~\ref{lem:abs} and get $1/n_0 \ll \phi \ll \e \ll \r$. 
In addition, we assume that $\phi \ll \e \ll \a$.
Suppose
\[
1/n_0 \ll \phi \ll \e \ll \a \ll \r \ll \mu, \beta \ll 1/D_F, 1/r, 1/k,
\]
and $n\ge n_0$ with $k \mid r n$.
We also require $n$ to be large enough such that we can apply Lemma~\ref{lem:keevash} with constants $\phi$, $2\r$ and $\beta/2$ in place of $\phi$, $\r$ and $\beta$.
Let $\cP$ and $J$ be as assumed in the theorem.
Our aim is to show that $J_k$ contains a perfect matching which $\alpha$-represents $F$.
We apply Lemma~\ref{lem:abs} with $\zeta = 1/k-\gamma$ and get a balanced absorbing set $W$ of size at most $\e n$. 
Let $V'=V(J)\setminus W$ and $n'=n-|W|/r\ge (1-\e)n$. Let $J'=J[V']$ and note that every cluster of $J'$ has $n'$ vertices. 
Since $\e<\r$, we have
\[
\delta^F(J')\ge \left( n', \left(\frac{k-1}{k} - 2\r\right) n', \left(\frac {k-2}k n - 2\r\right) n', \dots, \left(\frac 1k - 2\r \right) n'  \right).
\]
%
By applying Lemma~\ref{lem:keevash} with the constants chosen above on $J'$, we get that either $J'_k$ contains an $F$-balanced matching $M$ which covers all but at most $\phi n' < \phi n$ vertices of $V'$, or for some $p \in [k-1]$ and sets $S_i \subseteq V_i$ with $|S_i| = \lfloor pn'/k \rfloor$, we have $e(J'_{p+1}[S] ) \le \beta (n')^{p+1}/2$, where $S:=\bigcup_{i\in [r]}S_i$.
If the latter holds, then for each $i\in [r]$ we add arbitrary $\lfloor pn/k \rfloor - \lfloor pn'/k \rfloor \le \e n$ vertices from $V_i\setminus S_i$ to $S_i$, and denote the resulting set by $S_i'$.
Note that for each $i\in [r]$, $S_i'$ has exactly $\lfloor pn/k \rfloor$ vertices and 
\[
e(J_{p+1}[S'] )< \beta (n')^{p+1}/2 + r\e n \cdot n^p < \beta n^{p+1},
\]
by the choice of $\e$, where $S' = \bigcup_{i\in [r]}S_i'$.
This contradicts~\ref{item:i}.

Otherwise the former holds.
Let $U=V(J')\setminus V(M)$ and thus $|U|\le \phi n$ and $k\mid |U|$.
By the property of $W$, there is a perfect matching on $J_k[U\cup W]$.
The union of this perfect matching and $M$ gives a perfect matching $M'$ of $J_k$.
Moreover, 
for each $\bfi\in I(F)$, write $n_{\bfi}(M')$ as the number of edges in $M'$ with index vector $\bfi$ divided by the multiplicity of $\bfi$ in $I(F)$.
Then, since $M$ is an $F$-balanced matching that covers all but $|W| + |U|$ vertices, for any $\bfi\in I(F)$, we have
\[
n_{\bfi}(M') \ge \frac{r n - |W| - |U|}{k|F|} \ge \frac{r n}{k|F|} - \frac{ \e n + \phi n} k,
\]
as $|W|\le \e n$ and $|U|\le \phi n$.
On the other hand, this gives that
\[
n_{\bfi}(M') \le \frac{r n}{k|F|} + |F| \frac{ \e n + \phi n} k. 
\]
Thus $M'$ $\a$-represents $F$ because for different $\bfi, \bfi'\in I(F)$,
\[
\frac{n_{\bfi}(M')}{n_{\bfi'}(M')} \ge \frac{r - |F| (\e + \phi)}{r + |F|^2 (\e + \phi)} > 1-\a
\]
as $|F|\le D_F$ and $\e, \phi \ll \a \ll 1/D_F$.
The proof is completed.
\end{proof}

It remains to prove Lemmas~\ref{lem:keevash} and~\ref{lem:abs}.

\section{Proof of the Absorbing Lemma}

We prove Lemma~\ref{lem:abs} in this section. We first introduce some notation and auxiliary results.

\subsection{Tools}
We use the following notation introduced by Lo and Markstr\"om~\cite{LM1}.
Let $H$ be a $k$-graph on a vertex set $V$ with $|V|=n$.
Two vertices $u, v\in V$ are called \emph{$(\beta, i)$-reachable} in $H$ if there are at least $\beta n^{i k-1}$ $(i k-1)$-sets $S$ such that both $H[S\cup u]$ and $H[S\cup v]$ have perfect matchings. 
We say a vertex set $U$ is \emph{$(\beta, i)$-closed} in $H$ if any two vertices $u,v\in U$ are $(\beta, i)$-reachable in $H$.
Note that when we work with a given dense $k$-complex $J$, the reachability will be defined on $J_k$.

Let $H$ be a $k$-graph with a partition $\cP=\{V_1,\dots, V_d\}$.
For any $j\in [d]$ and $v\in V(H)$, let $\tilde{N}_{\beta, j}(v, H)$ be the set of vertices in the same part as $v$ that are $(\beta, j)$-reachable to $v$.
As usual, for a $k$-graph $H$ and a set $S\subseteq V(H)$, let $N_H(S)$ be the collection of $(k-|S|)$-sets $T$ in $V(H)\setminus S$ such that $S\cup T\in E(H)$; let $\deg_H(S):=|N_H(S)|$.
Given another set $W$, let $\deg_H(S, W) := |N_H(S)\cap \binom{W}{k-|S|}|$.

\begin{proposition}\label{prop:Nv}
Suppose that $1/n \ll \eta \ll \zeta, 1/r, 1/k$. 
Suppose $F$ is an allocation such that for every $j\in [r]$ there exists $f\in F$ such that $f^{-1}(j)\neq \emptyset$.
Let $\cP$ be a partition of a set $V$ into parts $V_1,\dots, V_r$ of size $n$, and $J$ be a $\cP F$-partite $k$-complex on $V$ satisfying $\delta^F(J)\ge (n, \zeta n, \dots, \zeta n)$.
Then for any $v\in V_j$, $j\in [r]$, $|\tilde{N}_{\eta, 1}(v, J_k)\cap V_j| \ge \delta_{k-1}^F(J) - \sqrt\eta n$.
\end{proposition}

\begin{proof}
Fix $j\in [r]$ and $f\in F$ such that $f^{-1}(j)\neq \emptyset$.
For any vertex $w\in V_j$, let $N_f(w)$ be the collection of $(k-1)$-sets $S\in N_{J_k}(w)$ such that $\bfi_\cP(S\cup w)=\bfi(f)$.
Pick a vertex $v\in V_j$.
Note that for any other vertex $u\in V_j$, $u\in \tilde{N}_{\eta, 1}(v, J_k)$ if $|N_{f}(u)\cap N_{f}(v)|\ge \eta (r n)^{k-1}$.
By double counting, we have
\[
 \sum_{S\in N_f(v)} \deg_{J_k}(S, V_j) < |\tilde{N}_{\eta, 1}(v, J_k)|\cdot |N_f(v)|+n\cdot \eta {(r n)}^{k-1}.
\]
For any $S\in N_f(v)$ in the above inequality, since $S$ can be constructed by following a certain permutation of $f$, it holds that $\deg_{J_k}(S, V_j) \ge \delta_{k-1}^F(J)$. Moreover, we have that
\[
|N_f(v)| \ge \frac{1}{(k-1)!} (\zeta n)^{k-1}\ge \sqrt\eta (r n)^{k-1}.
\]
Thus, $|\tilde{N}_{\eta, 1}(v, J_k)|> \delta_{k-1}^F(J) - {\eta r^{k-1} n^k}/{|N_f(v)|}\ge \delta_{k-1}^F(J) - \sqrt\eta n$ as desired.
\end{proof}

We use the following lemma from \cite[Lemma 6.3]{HT}, which is similar to a result first appeared in~\cite{Han14_poly}. 


\begin{lemma}\cite{HT}\label{lem:P}
Given $\delta>0$, integers $k \ge 2$ and $0<\a \ll \delta, 1/k$, there exists a constant $\beta>0$ such that the following holds for all sufficiently large $n$. 
Assume $H$ is an $n$-vertex $k$-graph  and $S\subseteq V(H)$ is such that $|\tilde{N}_{\a, 1}(v,H)\cap S| \ge \delta n$ for any $v\in S$. 
Then we can find a partition $\cP$ of $S$ into $V_1,\dots, V_r$ with $r\le 1/\delta$ such that for any $i\in [r]$, $|V_i|\ge (\delta - \a) n$ and $V_i$ is $(\beta, 2^{\lfloor 1/\delta \rfloor-1})$-closed in $H$.
\end{lemma}

%


Fix an integer $i>0$ and let $H$ be a $k$-graph.
For a $k$-vertex set $S$, we call a set $T$ an \emph{absorbing $i$-set for $S$} if $|T|=i$ and both $H[T]$ and $H[T\cup S]$ contain perfect matchings. 
We use the absorbing lemma from \cite[Lemma 3.4]{Han14_poly} with some quantitative changes, which follows from the original proof (a similar formulation appears in \cite{Han15_mat}).


\begin{lemma}\cite{Han14_poly}\label{lem:abs0}
Suppose $1/n \ll \phi \ll \beta, \mu \ll 1/k, 1/r, 1/t$
and $H$ is a $k$-graph on $n$ vertices.
Suppose $\cP=\{V_1, \dots, V_r\}$ is a partition of $V(H)$ such that for $i\in [r]$, $V_i$ is $(\beta, t)$-closed. Then there is a family $\F_{abs}$ of disjoint $tk^2$-sets with size at most $\beta n$ such that $H[V(\F_{abs})]$ contains a perfect matching and every $k$-vertex set $S$ with $\bfi_{\cP}(S)\in I_{\cP}^{\mu} (H)$ has at least $\phi n$ absorbing $t k^2$-sets in $\F_{abs}$. 
\end{lemma}

Let $I$ be a set of $k$-vectors in $\mathbb Z^r$ and let $L$ be the lattice generated by $I$.
Then for any $k$-vector $\bfw\in L$, there exist integers $a_\bfv$ for each $\bfv\in I$ such that
\[
\bfw = \sum_{\bfv\in I} a_\bfv \bfv.
\]
Let $C(k, r, I)$ be the maximum of $\max_{\bfv\in I}|a_\bfv|$ over all $r^k$ choices for $k$-vectors $\bfw$.
Then let $C:=C(k,r)$ be the maximum of $C(k, r, I)$ over choices for $I$.

Now we are ready to prove Lemma~\ref{lem:abs}.
Here is an outline of the proof.
Let $J$ be a $\cP F$-partite complex satisfying the assumptions of the lemma.
We use Lemma~\ref{lem:P} on $J_k$ to build a partition $\cP=\{V_1,\dots, V_d\}$ of $V(J)$ for some constant $d$ such that every $V_i$ is $(\beta, t)$-closed for some $\beta>0$ and integer $t\ge 1$. 
Since the lattice is complete by~\ref{item:i}, Lemma~\ref{lem:abs0} gives the desired absorbing set.


\begin{proof}[Proof of Lemma~\ref{lem:abs}]
Take additional constants $\beta, \eta, C>0$ such that $\phi \ll \beta \ll \epsilon, \eta \ll \zeta \ll 1/C$, where $C:=C(k,r)$ is defined as above.
We first note that because $F$ is $(k,r)$-uniform, the assumption of Proposition~\ref{prop:Nv} is satisfied.
Thus, by Proposition~\ref{prop:Nv}, for any $v\in V_j$, $j\in [r]$, $|\tilde{N}_{\eta, 1}(v, J_k)\cap V_j| \ge \delta_{k-1}^F(J) - \sqrt\eta n \ge (\zeta - \sqrt \eta)n$.
Therefore we can apply Lemma~\ref{lem:P} on $J_k$ with $S=V_i$ and $\delta=(\zeta-\sqrt\eta)/r$ (note that $|V(J)|=r n$) for $i\in [r]$ respectively.
This gives a partition
\[
\cP'=\{V_{11}, V_{12}, \dots, V_{1a_1}, V_{21}, V_{22}, \dots, V_{2a_2},\dots, V_{r1},\dots, V_{ra_r}\}
\]
such that each $V_{ij}$ is $(\beta, 2^{\lfloor 1/\delta \rfloor-1})$-closed in $J_k$, $|V_{ij}|\ge (\delta - \eta)r n\ge (\zeta-\mu)n$ and $V_i=\bigcup_{j\in [a_i]} V_{ij}$.
Note that the lower bound of $|V_{ij}|$ implies $a_i\le 1/(\zeta - \mu)$ and clearly $\cP'$ refines $\cP$.
Write $I:=I_{\cP'}^\mu(J_k)$ and $L:=L_{\cP'}^\mu(J_k)$.
Thus by the assumption that Theorem~\ref{thm:711}~\ref{item:iii} holds, $L$ is complete.

We apply Lemma~\ref{lem:abs0} on $J_k$ with $(Cr^k+1)\phi$ in place of $\phi$ and get the family $\F_{abs}$.
Denote the perfect matching on $V(\F_{abs})$ by $M_1$.
Take a matching $M_2:=\bigcup_{\bfi\in I}M_\bfi$, where each $M_\bfi$ is a matching of $C\phi n$ edges all of index vector $\bfi$.
Note that we can greedily construct $M_2$ because $C\phi |I|n\le r^kC\phi n< \mu n$ and there are at least $\mu n^k$ edges with allocation function $f$ for each $f\in F$.
Since $J$ is $\cP F$-partite, every edge of $J_k$ corresponds to some allocation function $f\in F$.
We extend $M_1\cup M_2$ greedily to an $F$-balanced matching.
Since $|F|\le D_F$, the resulting matching would have size at most $D_F(|M_1|+|M_2|)$, so it can be constructed greedily as $|M_1|+|M_2|\le tk \beta n + r^k C\phi n < \mu n/D_F$ and there are at least $\mu n^k$ edges with allocation function $f$ for each $f\in F$.
Let $W$ be the vertex set of the resulting matching.
It holds that $|W|\le kD_F(|M_1|+|M_2|) < \epsilon n$ and it remains to verify the absorption property.

Let $U \subseteq V(J) \setminus W$ be any set such that $|U|\le \phi n$ and $k\mid |U|$.
We arbitrarily partition $U$ into $k$-sets $S_1,\dots, S_t$, where $t:=|U|/k$.
Since $L$ is complete, for each $i\in [t]$, there exists $\{a_\bfv\in \mathbb Z: \bfv\in I\}$ such that
\[
\bfi(S_i) = \sum_{\bfv\in I} a_\bfv \bfv
\]
and $|a_\bfv|\le C$, $\bfv\in I$.
For each $\bfv\in I$, one can rewrite $a_\bfv = b_\bfv - c_\bfv$ such that one of $b_\bfv$ and $c_\bfv$ is in $[C]$ and the other is zero.
Thus we obtain
\begin{equation}\label{eq:bc}
\bfi(S_i) + \sum_{\bfv\in I} c_\bfv \bfv = \sum_{\bfv\in I} b_\bfv \bfv.
\end{equation}
That is, we can take $c_\bfv$ edges of index vector $\bfv$ for each $\bfv\in I$ from $M_2$, and decompose the union of these edges and $S_i$ as a collection of $k$-sets consisting of $b_\bfv$ $k$-sets for each $\bfv\in I$.
Repeating this for all $i\in [t]$ (by adding disjoint edges from $M_2$), we obtain a collection $\mathcal T$ of disjoint $k$-sets with index vector in $I$.
This is possible because it consumes at most $C\phi n$ edges for each $\bfi\in I$ from $M_\bfi\subseteq M_2$; moreover, $|\mathcal T|\le (1+C|I|)\phi n\le (1+Cr^k)\phi n$.
Thus, these $k$-sets in $\mathcal T$ can be absorbed by members of $\F_{abs}$.
This shows that $J_k[U\cup W]$ has a perfect matching and we are done.
\end{proof}

\section{A reproof of Lemma~\ref{lem:keevash}}

In this section we give another proof of Lemma~\ref{lem:keevash}, which avoids the use of hypergraph regularity.

\subsection{Perfect fractional matching}

A function $g\colon E(H) \to[0,1]$ is called a \emph{perfect fractional matching}
if $\sum_{e\in H\colon v\in e} g(e)= 1$ for every $v\in V(H)$.
Let $F$ be a $(k,r)$-uniform allocation.
Given a $\cP F$-partite $k$-system $J$, we say a perfect fractional matching in $J$ is \emph{$F$-balanced}, if $\sum_{e\in H: \bfi(e)=\bfi} g(e)/m_{\bfi}$ is a constant over all $\bfi\in I(F)$, where $m_{\bfi}$ is the multiplicity of $\bfi$ in $I(F)$.

We use the following result on perfect fractional matchings from~\cite[Lemma 7.2]{KM1}, which says that the assumptions in~Lemma~\ref{lem:keevash} guarantees an $F$-balanced perfect fractional matching.
It is proved by a careful application of the well-known Farkas' Lemma on the solvability of linear inequalities.

\begin{lemma}\cite{KM1} \label{lem:frac_mat}
Suppose that $1/n \ll \r \ll \beta, 1/D_F, 1/r, 1/k$.
Suppose $F$ is a $(k,r)$-uniform connected allocation with $|F|\le D_F$.
Let $\cP$ be a partition of a set $V$ into parts $V_1,\dots, V_r$ of size $n$, and $J$ be a $\cP F$-partite $k$-system on $V$ 
satisfying \eqref{eq:deg3} and~\ref{item:i} in Theorem~\ref{thm:711}.
Then $J_k$ contains an $F$-balanced perfect fractional matching $g$.
\end{lemma}

Our first goal is to prove the following result, which says that under the same assumption, one can find many `weight-disjoint' (see~\eqref{eq:multi}) perfect fractional matchings.
These perfect fractional matchings will enable us to choose a random subgraph $H$ of $J_k$, which is almost regular.
Moreover,~\eqref{eq:multi} implies that the maximum pair degree of $H$ is small.
It is known that such a $k$-graph contains an almost perfect matching (see Theorem~\ref{thm:Pippenger} below).
However, such a matching may not be {$F$-balanced}.
To overcome this issue, we define another auxiliary hypergraph $H''$, which also satisfies the assumptions of Theorem~\ref{thm:Pippenger} and more importantly, any matching of $H''$ can be decomposed into an $F$-balanced matching of $J_k$.

\begin{lemma} \label{lem:multi_frac}
Suppose that $1/n \ll \r \ll \beta, 1/D_F, 1/r, 1/k$.
Suppose $F$ is a $(k,r)$-uniform connected allocation with $|F|\le D_F$.
Let $\cP$ be a partition of a set $V$ into parts $V_1,\dots, V_r$ of size $n$, and $J$ be a $\cP F$-partite $k$-system on $V$ 
satisfying \eqref{eq:deg3} and~\ref{item:i} in Theorem~\ref{thm:711}.
Then $J_k$ contains $\ell=\r n$ $F$-balanced perfect fractional matchings $g_1,\dots, g_\ell$ such that
 \begin{equation}
    \label{eq:multi}
    \sum_{e: \,u, v\in e}\sum_{i=1}^\ell g_i(e)\le 2
    \text{ for every pair }uv\in \binom V2.
  \end{equation}
\end{lemma}


\begin{proof}
We define a weight function $w$ on the pairs of $V$ to track the available weight on edges.
At the beginning, set $w(e):=2$ for all $e\in \binom V2$.
Throughout the process, let $G\subseteq \binom V2$ be the set of pairs $e$  such that $w(e)\ge 1$.
Let $J'\subseteq J$ be the $k$-system with edges not supported on $G$ removed. 
We shall iteratively apply Lemma~\ref{lem:frac_mat} to $J'$ to find $\ell=\r n$
$F$-balanced perfect fractional matchings $g_1,\dots, g_\ell$. In doing so we will
iteratively update the weights of the pairs in $\binom V2$, that is, for each $uv\in \binom V2$, we let $w(uv):=w(uv) - \sum_{e\in J_k: u,v\in e}g_i(e)$.

Consider any intermediate step $\ell'\le \ell$.
Note that any pair $uv$ which is removed during the process is because $w(uv)<1$, namely, at that point the weight we have chosen is at least 1.
Therefore, as for any $v$, the weight we have chosen so far is $\sum_{i\in [\ell']}\sum_{e\in H\colon v\in e} g_i(e)= \ell'\le \ell$, there are at most $\ell$ pairs containing $v$ which are not in $G$,
that is, $\delta^F_1(J') \ge \delta^F_1(J) - \ell$.
Moreover, for each $3\le i\le k$ and $R\in J_{i-1}'$, the number of edges $T\in J_i$ such that $R\subseteq T$ and $w(e)<1$ for some pair $e\subset T$ is at most $(i-1)\ell$ (because there are $i-1$ pairs $e$ to be considered), that is, $\delta^F_{i-1}(J') \ge \delta^F_{i-1}(J) - (i-1)\ell$.
Since $\ell=\r n$, we obtain
\[
\delta^F(J')\ge \left( n, \left(\frac{k-1}{k} - \r\right) n-\r n, \left(\frac {k-2}k - \r\right) n-2\r n, \dots, \left(\frac 1k -\r \right) n -(k-1)\r n  \right)
\]
and for each $p\in [k-1]$, $e(J_{p+1}\setminus J'_{p+1})\le p\r n\cdot n^p<\beta n^{p+1}/2$.
Moreover, suppose there exists $p\in [k-1]$ and sets $S_i\subseteq V_i$ with $S_i = \lfloor pn/k \rfloor$ for all $i\in [r]$, we have $e(J'_{p+1}[S] )\le \beta n^{p+1}/2$, where $S:=\bigcup_{i\in [r]}S_i$.
Then we have $e(J_{p+1}[S])\le e(J'_{p+1}[S] )+ e(J_{p+1}\setminus J'_{p+1})\le \beta n^{p+1}$, contradicting our assumptions.
So we can apply Lemma~\ref{lem:frac_mat} with $k\r$ in place of $\r$ and $\beta/2$ in place of $\beta$ to find an $F$-balanced perfect fractional matching in the current weighted $k$-system.  
Note that~\eqref{eq:multi} holds because every $w(uv)\ge 0$ at the end of the process.
\end{proof}

We will use the following theorem of Frankl and R\"odl~\cite{FrRo}
(see also R\"odl~\cite{VR85} and Alon and
Spencer~\cite[Theorem~4.7.1]{AS16}), which asserts the existence of an
almost perfect matching in `pseudorandom' hypergraphs.
For a hypergraph $H$, let $\Delta_2(H): = \max_{u,v\in V(H)}\deg(u v)$.
For reals $a, b$ and $c$, we write $a=(1\pm b)c$ for $(1-b)c\le a\le (1+b)c$.
Theorem~\ref{thm:Pippenger} has been applied to similar problems, see our concluding remarks.

\begin{theorem}\cite{FrRo}\label{thm:Pippenger}
  For every integer $k\ge 2$ and a real $\e>0$, there exist $\tau$ and $d_0$ such that
  for every $n\ge D\ge d_0$ the following holds.
  Every $k$-uniform hypergraph on a set $V$ of $n$ vertices which satisfies the following conditions:
  \begin{enumerate}
  \item for all vertices $v\in V$, we have $\deg_H(v)= (1\pm \tau)D$ and
  \item $\Delta_2(H)\le \tau D$
  \end{enumerate}
  contains a matching covering all but at most $\e n$ vertices. 
\end{theorem}

Combining Lemma~\ref{lem:multi_frac} and the following lemma we obtain a reproof of Lemma~\ref{lem:keevash}.

\begin{lemma} \label{lem:frac_int}
Suppose that $1/n \ll \phi \ll \r \ll \beta, 1/D_F, 1/r, 1/k$.
Suppose $F$ is a $(k,r)$-uniform connected allocation with $|F|\le D_F$.
Let $\cP$ be a partition of a set $V$ into parts $V_1,\dots, V_r$ of size $n$, and $J$ be a $\cP F$-partite $k$-system on $V$.
Suppose $J_k$ contains $\ell=\r n$ $F$-balanced perfect fractional matchings $g_1,\dots, g_\ell$ such that~\eqref{eq:multi} holds.
Then $J_k$ contains an $F$-balanced matching $M$ which covers all but at most $\phi n $ vertices of~$V$.
\end{lemma}

\begin{proof}
Write $d=k|F|\le kD_F$.
Let $\tau_0$ be returned from Theorem~\ref{thm:Pippenger} with $\e=\phi/r$ and $d$ in place of $k$.
Let $\tau=\tau_0/(2d)$.
For each edge $e\in J_k$, define $g(e) := \sum_{i\in [\ell]} g_i(e)/2$.
In particular, $g(e)\le 1$ for any $e\in J_k$.
Moreover, since each $g_i$ is $F$-balanced, for every $\bfi\in I(F)$, the number
$\sum_{e:\, \bfi(e)=\bfi} g_i(e)/m_{\bfi} = (r n/k)/|F|= r n/d$ is a constant, where $m_{\bfi}$ is the multiplicity of $\bfi$ in $I(F)$.
Thus, $\sum_{e:\, \bfi(e)=\bfi} g(e)/m_{\bfi} = \ell r n/d$.

Next we select a random subgraph $H$ of $J_k$, by including each edge $e\in J_k$ independently with probability $g(e)$.
Note that for every $v\in V(J)$, $\sum_{e: v\in e}g(e)=\ell$.
So we have $\mathbb E[\deg_H(v)] =\sum_{e: v\in e}g(e)= \ell$ for every $v\in V(J)$ and by~\eqref{eq:multi}
\[
\mathbb E[\deg_H(uv)] = \sum_{e: \,u, v\in e} g(e)= \sum_{e: \,u, v\in e}\sum_{i=1}^\ell g_i(e)/2 \le 1
\]
for any $u,v \in V(J)$. 
Let $Y_{\bfi}$ be the random variable that counts the number of edges in $H$ with index vector $\bfi\in I(F)$.
Then $\mathbb EY_{\bfi}=\sum_{e:\, \bfi(e)=\bfi} g(e) = m_{\bfi}\ell r n/d$.
By standard concentration results (e.g., Chernoff's bound), we infer that with positive probability,
\begin{enumerate}[label=(\alph*)]
\item $\deg_H(v)  = (1\pm \tau) \ell$ for every $v\in V(J)$,  \label{item:a}
\item $\Delta_2(H)\le 3\log n$, and \label{item:b}
\item the number of edges of $H$ with index vector $\bfi$ is $(1\pm \tau/2)m_{\bfi} \ell r n/ d$ for every $\bfi\in I(F)$. \label{item:c}
\end{enumerate}
So there is a $k$-graph $H$ satisfying all these properties.

For each $\bfi\in I(F)$, we split the edges of $H$ with index vector $\bfi$ into $m_{\bfi}$ color classes as equal as possible arbitrarily.
This defines a coloring $E_1, \dots, E_{|F|}$ of the edges of $H$ by $|F|$ colors, each with $(1\pm \tau/2) \ell r n/ d$ edges (by~\ref{item:c}).
Let $H'$ be a $d$-graph on $V(J)$ such that the edges of $H'$ are the $d$-sets that are the union of $|F|$ disjoint edges $e_1,\dots, e_{|F|}$ in $H$, one from each of the $|F|$ color classes.
Thus, for every $v\in V(J)$, by~\ref{item:a}, we have $\sum_{i\in [|F|]} \deg_{E_i}(v) = \deg_H(v) = (1\pm \tau) \ell$.
Therefore we obtain
\[
\deg_{H'}(v) = \sum_{i\in [|F|]} \deg_{E_i}(v) \left((1\pm \tau) \ell \frac{r n} d\right)^{d-1} = (1\pm \tau) \ell \left((1\pm \tau) \ell \frac{r n} d\right)^{d-1} = (1\pm d\tau) \ell \left(\ell \frac{ r n} d\right)^{d-1}
\]
and $\Delta_2(H')\le 4(\log n) \left({\ell r n} /d\right)^{d-1}$ by~\ref{item:b}.

Finally, we select a subgraph $H''\subseteq H'$ by including each edge of $H'$ independently with probability $(\ell rn/d)^{1-d}$.
Again, by standard concentration results, we infer that with positive probability,
\begin{itemize}
\item $\deg_{H''}(v)  = (1\pm 2d\tau) \ell = (1\pm 2d\tau)\r n$ for every $v\in V(J)$, and
\item $\Delta_2(H'')\le 12(\log n)$.
\end{itemize}
So we can apply Lemma~\ref{thm:Pippenger} on $H''$ with $\e = \phi/r$ and obtain a matching $M'$ covering all but at most $\e rn = \phi n$ vertices.
Recall that as $F$ is $(k,r)$-uniform, for every $i\in [k]$ and $j\in [r]$ there are $|F|/r$ functions $f\in F$ with $f(i)=j$.
Thus, any edge of $M'$ contains exactly $d/r$ vertices from each $V_i$, $i\in [r]$ and can be decomposed into an $F$-balanced matching of size $|F|$ in $J_k$.
Let $M$ be the matching in $J_k$ given by $M'$.
By our construction, $M$ is $F$-balanced.
\end{proof}

\section{Concluding Remarks}

In this note we give an alternative proof of the main technical result~\cite[Theorem 7.1]{KM1}, which avoids the hypergraph blow-up lemma and the hypergraph regularity method. 
This allows us to obtain regularity-free proofs of the results in~\cite[Section 2]{KM1}, as well as the subsequent applications in~\cite{HLTZ_K4, KKM13,KM1,KeMy2}.
We did draw substantial notation and ideas (e.g. the theory of edge-lattice) from~\cite{KM1} as well as in~\cite{Han14_poly}.
The proof we present is not self-contained, but all the proofs of the three lemmas we cited in Section 4 are one to two pages long; the only substantial one is Lemma~\ref{lem:frac_mat} (\cite[Lemma 7.2]{KM1}), whose proof is six-page long. 

The authors of~\cite{KM1} also exploited their method and gave a general result which almost entirely dispenses with degree assumptions, assuming that the reduced $k$-complex, after applying the regularity lemma, has a perfect matching (see~\cite[Chapter 9]{KM1} for more details).
Here we present a similar result given by our method.
Since we tend to avoid the regularity method, we formulate the result as follows, assuming the existence of perfect fractional matchings in the almost spanning subcomplex.

\begin{theorem}\label{thm:general}
Let $1/n \ll \r ,\e , \mu \ll \zeta, 1/k$ and $k\mid n$.
Let $J$ be a $k$-complex on $n$ vertices satisfying $\delta(J)\ge (n, \zeta n, \dots, \zeta n)$ and
\begin{enumerate}[label=$(\roman*)$]
\item every induced subcomplex $J'$ of $J$ on at least $n-\e n$ vertices has at least $\ell=\gamma |V(J')|$ perfect fractional matchings $g_1,\dots, g_{\ell}$ satisfying~\eqref{eq:multi}. 
\item $L_{\cP}^{\mu}(J_k)$ is complete for any partition $\cP$ of $V(J)$ whose parts each have size at least $(\zeta-\mu) n$.
\end{enumerate}
Then $J_k$ contains a perfect matching.
\end{theorem}

\begin{proof}
Choose an additional constant $1/n \ll \phi \ll \r$.
Apply Lemma~\ref{lem:abs} to obtain an absorbing set $W$ of size $\e n$.
Let $J'$ be the induced subcomplex of $J$ on $V(J)\setminus W$ and $J'$ contains those perfect fractional matchings by $(i)$.
Then Lemma~\ref{lem:frac_int} gives a matching that leaves a set $U$ of vertices in $J'$ uncovered, where $k\mid |U|$ and $|U|\le \phi n$.
Absorbing these vertices by $W$ gives a perfect matching in $J_k$.
\end{proof}

It is also possible to replace the minimum degree sequence condition in Theorem~\ref{thm:general} above by ``every vertex is in at least $\zeta n^{k-1}$ edges in $J$'' by applying a `shaving' lemma in~\cite[Lemma 3.3]{HZ18}, which finds a spanning subcomplex $J'$ and a set $Y$ with $|Y|\le \zeta^2 n$ such that every $(k-1)$-set in $V(J)\setminus Y$ has degree either $0$ or at least $\eta n$ in $J'$ for some $\eta>0$ and every vertex in $V(J)\setminus Y$ is still in at least $(\zeta /2) n^{k-1}$ edges in $J'$.
Then following the proof of Lemma~\ref{lem:abs}, we can build the absorbing set $W$ which can absorb $k$-sets in $V(J)\setminus Y$.
We add to the absorbing set a small matching that covers $Y\setminus W$ (which is possible as $|Y|<\zeta n$), and the resulting set will have the absorbing property as in Lemma~\ref{lem:abs}.
We omit the details.


Our use of Theorem~\ref{thm:Pippenger} is not novel.
In fact, Alon et~al.~\cite{AFHRRS} also used Theorem~\ref{thm:Pippenger} to turn a perfect fractional matching into an almost perfect matching via probabilistic arguments.
The main difference is that their minimum degree condition is slightly above that of the space barriers; this slight difference enables them to carry out the argument without involving approximate space barriers.
In particular, they applied the fractional matching result on many sets of order around $n^{0.2}$.
In contrast, in our case, we may be in the situation that such a set of order $n^{0.2}$ is close to the space barriers.
Then it is not clear how to argue that this would happen only if the \emph{host} graph is close to the space barriers (if the small set has order $\Omega(n)$, it is straightforward to argue this by concentration results).
Instead, we use the new argument developed in~\cite{HKP}, which goes around this issue by finding `weight-disjoint' perfect fractional matchings.
Since each fractional matching is found in the original hypergraph, we can easily deal with the case when space barriers arise.

At last, unlike the regularity method, the probabilistic method usually gives moderate bounds on the order of the hypergraph.
In fact, the regularity method requires the order of the hypergraph to be at least $n_0$ which is a tower-type function (in fact, the $k$-th Ackermann function) of the regularity parameter $1/\e$.
To make the regularity useful one must have $\e < \r$, the main constant in the degree condition (here $k$ is usually considered as a constant).
In contrast, our probabilistic arguments (mainly by Chernoff's bound) only requires $n_0$ to be polynomial in $1/\r$.
For example, in the proof of Lemma~\ref{lem:frac_int}, when we choose a random subgraph of $J_k$ by including edges uniformly at random, the failure probability for each event will be of form $e^{-\r^c n}$ for some constant $c$.
Then if we require $n^2$ events to hold simultaneously, this will require $n^2 e^{-\r^c n} < 1$, in particular, $\r^c n \ge \sqrt n$ will be sufficient.


\section*{Acknowledgement}
We thank Yoshiharu Kohayakawa, Richard Mycroft and Yi Zhao for helpful discussions.
We also thank two referees for comments that improve the presentation of this paper.

\bibliographystyle{abbrv}
\bibliography{Bibref}

\end{document}